\documentclass[draft]{jotart}

\usepackage{amssymb,amsmath}
\usepackage{bm}
\usepackage[all,cmtip]{xy}
\usepackage{amscd}

\usepackage{ifthen}
\usepackage{graphicx}
\usepackage{color}

\newcommand{\bea}{\begin{eqnarray}}
\newcommand{\eea}{\end{eqnarray}}

\newcommand{\vp}{\varphi}

\newcommand{\D}{\mathbb{D}}

\def\textmatrix#1&#2\\#3&#4\\{\bigl({#1 \atop #3}\ {#2 \atop #4}\bigr)}
\def\dispmatrix#1&#2\\#3&#4\\{\left({#1 \atop #3}\ {#2 \atop #4}\right)}
\newcommand{\be}{\begin{equation}}
\newcommand{\ee}{\end{equation}}
\newcommand{\ben}{\begin{eqnarray*}}
\newcommand{\een}{\end{eqnarray*}}

\newcommand{\NI}{\noindent}

\newcommand{\bi}{\begin{itemize}}
\newcommand{\ei}{\end{itemize}}

\newcommand{\es}{{\mathcal S}}
\newcommand{\OO}{{\mathcal O}}
\newcommand{\ds}{\displaystyle}

\theoremstyle{proclaim}

\theoremstyle{statement}

\theoremstyle{definition}

\theoremstyle{plain}

\newtheorem{thm}{Theorem}[section]
\newtheorem{cor}[thm]{Corollary}

\newtheorem{prop}[thm]{Proposition}
\theoremstyle{definition}

\newtheorem{rem}[thm]{Remark}

\numberwithin{equation}{section}

\let\phi=\varphi

\begin{document}

\issueinfo{00}{0}{0000}
\commby{F.-H. Vasilescu}
\pagespan{000}{000}
\date{May 15, 2020}

\title[Invariant subspaces of composition operators]{Beurling type invariant subspaces of composition operators}


\author[S. Bose, P. Muthukumar{\protect \and} J. Sarkar]{Snehasish Bose, P. Muthukumar {\protect \and} Jaydeb Sarkar}

\address{Indian Statistical Institute, Statistics and Mathematics Unit, 8th Mile, Mysore Road, Bangalore, 560059, India}
\email{bosesonay@gmail.com}

\address{Indian Statistical Institute, Statistics and Mathematics Unit, 8th Mile, Mysore Road, Bangalore, 560059, India}
\email{pmuthumaths@gmail.com}

\address{Indian Statistical Institute, Statistics and Mathematics Unit, 8th Mile, Mysore Road, Bangalore, 560059, India}
\email{jay@isibang.ac.in, jaydeb@gmail.com}




\begin{abstract}
Let $\mathbb{D}$ be the open unit disk in $\mathbb{C}$, let $H^2$ denote the Hardy space on $\mathbb{D}$ and let $\varphi : \mathbb{D} \rightarrow \mathbb{D}$ be a holomorphic self map of $\mathbb{D}$. The composition operator $C_{\varphi}$ on $H^2$ is defined by
\[
(C_{\varphi} f)(z)=f(\varphi(z)) \quad \quad (f \in H^2,\, z \in \mathbb{D}).
\]
Denote by $\mathcal{S}(\mathbb{D})$ the set of all functions that are holomorphic and bounded by one in modulus on $\mathbb{D}$, that is
\[
\mathcal{S}(\mathbb{D}) = \{\psi \in H^\infty(\mathbb{D}): \|\psi\|_{\infty} := \sup_{z \in \mathbb{D}} |\psi(z)| \leq 1\}.
\]
The elements of $\mathcal{S}(\mathbb{D})$ are called Schur functions. The aim of this paper is to answer the following question concerning invariant subspaces of composition operators: Characterize $\varphi$, holomorphic self maps of $\mathbb{D}$, and inner functions $\theta \in H^\infty(\mathbb{D})$ such that the Beurling type invariant subspace $\theta H^2$ is an invariant subspace for $C_{\varphi}$. We prove the following result: $C_{\varphi} (\theta H^2) \subseteq \theta H^2$
if and only if
\[
\frac{\theta \circ \varphi}{\theta} \in \mathcal{S}(\mathbb{D}).
\]
This classification also allows us to recover or improve some known results on Beurling type invariant subspaces of composition operators.
\end{abstract}

\begin{subjclass}
Primary: 47B33; Secondary: 30H10, 47B38, 30D55, 46E15\\
\phantom{MSC20102010} 47A15, 46E22
\end{subjclass}
\begin{keywords}
Composition operators, invariant subspaces, inner functions, Blaschke products, Schur functions, singular inner functions, Hardy space
\end{keywords}

\maketitle

\section{Introduction}

The invariant subspace problem \cite{HH}, one of the most important open problems in linear analysis, asks if every bounded linear operator on a separable Hilbert space has a non-trivial closed invariant subspace. This problem has an equivalent form which turns it into a more concrete function theoretic problem. To be more specific, let $\D$ be the open unit disk in $\mathbb{C}$, let $H^2$ denote the Hardy space on $\mathbb{D}$ and let $\varphi : \mathbb{D} \rightarrow \mathbb{D}$ be a holomorphic self map of $\D$. The composition operator $C_{\varphi}$ on $H^2$ is defined by $C_{\vp} f = f \circ \vp$, that is
\[
(C_{\varphi} f)(z)=f(\varphi(z)),
\]
for all $f \in H^2$ and $z \in \mathbb{D}$. Littlewood's subordination principle \cite{Shapiro:Book}  implies that $C_{\vp}$ is a bounded operator on $H^2$ and
\[
\|C_\vp\| \leq \sqrt{\frac{1 + |\vp(0)|}{1 - |\vp(0)|}}.
\]
By \cite{NRW,ISP-CO}, the invariant subspace problem for Hilbert space operators can be reformulated by considering any, fixed hyperbolic disc automorphism $\vp$. Indeed, in the aforementioned papers it is shown that the (unknown) fact that any Hilbert
space operator acting on a complex, infinite-dimensional, separable space always has proper invariant subspaces, is equivalent to the fact that the only minimal invariant subspaces of $C_{\vp}$ are the $1$-dimensional eigenspaces.

While descriptions of invariant subspace lattices of composition operators exist, (see for instance \cite{MPS}), the result in \cite{NRW,ISP-CO} referred above implies that an automorphic, hyperbolic composition operator has a very rich and complicated invariant subspace lattice, and so, one way to understand it, would be by describing sublattices, for instance that consisting of joint invariant subspaces of $C_{\vp}$ and $M_z$ on $H^2$. Here $M_z$ denote the unilateral forward shift operator or the multiplication operator induced by the coordinate function $z$ on $H^2$. The closed invariant subspaces of $M_z$ are called \textit{Beurling type subspaces} (or \textit{Beurling subspaces}).

The initiative of studying joint invariant subspaces of $C_{\vp}$ and $M_z$ is recent \cite{Inv2005}. Among the papers bringing up convincing arguments that a line of research like that is interesting we specify \cite{Inv2014}, \cite{Inv2001} and \cite{Inv2015}. In those papers, it is observed and proved, that classical theorems in function theory, most notoriously, the Julia-Carath\'eodory theorem, can be understood in terms of the action of composition operators on Beurling subspaces. That theorem addresses the existence of angular derivatives in the sense of Constantin
Carath\'eodory, and the authors of \cite{Inv2014} and \cite{Inv2015} observe that the existence of such an angular derivative of some analytic self map $\vp$ of the unit disc is equivalent to the fact that $C_{\vp}$ maps certain Beurling subspaces induced by some atomic singular inner functions into similar (not necessarily identical) spaces.

It is evident now that the joint invariant subspace problem of $C_{\vp}$ and $M_z$ introduces also a lot of additional structure of holomorphic self maps and inner functions. Indeed, the notion of inner functions arose as a result of the representations of shift invariant subspaces of the Hardy space. Recall that an \text{inner function} is a function $\theta \in H^2$ whose radial limits have modulus one a.e. on $\partial \D$. A classical result of A. Beurling \cite{InvSub-MO} classifies the invariant subspaces of $M_z$ as follows:

\noindent \textsf{Beurling's Theorem:} Let $\es \neq \{0\}$ be a closed subspace of $H^2$. Then $\es$ is invariant under $M_z$ if and only if there exists an inner function $\theta$ (unique up to a scalar factor of unit modulus) such that
\[
\es = \theta H^2.
\]

Among many other results, Matache \cite{Inv2015} proved that for every holomorphic self map $\vp$ of $\D$ there exists a non-trivial $M_z$-invariant closed subspace $\mathcal{S} \subsetneqq H^2$ (depending on $\vp$) such that $C_{\vp} \mathcal{S} \subseteq \mathcal{S}$ (also see Theorem \ref{thm-nontrivial inv sub} for a new proof).

At the present stage, it is also worthwhile to recall the following problem \cite[Problem 1]{Inv2015}: \textit{Characterize in measure theoretical terms when $C_{\vp}(\theta_1 H^2) \subseteq \theta_2H^2$, where $\theta_1$ and $\theta_2$ are singular inner functions}.
If $\theta_1=\theta_2$, this problem becomes an invariant subspace problem, namely, \textit{``When is
a Beurling subspace induced by a singular inner function left invariant by a composition operator''?} We point out that \cite[Problem 1]{Inv2015} is solved in \cite[Corollary 2.15]{Inv2015} in the particular case of singular inner functions induced by purely atomic measures.

Typical results and proofs in this direction (including the ones mentioned above) often involves analytic properties of $\vp$ like (Denjoy-Wolff) fixed points and derivative of $\vp$ at fixed points. However, due to the complex classificational structure of (bi-)holomorphic self maps of $\D$, most known results are case-specific. But, from a more general point of view, we prove the following result: Let $\varphi$ be a holomorphic self map of $\mathbb{D}$, and let $\theta \in H^\infty(\mathbb{D})$ be an inner function. Then, the Beurling type invariant subspace $\theta H^2$ is invariant under $C_{\varphi}$ (that is, $C_{\vp} (\theta H^2) \subseteq \theta H^2$) if and only if
\[
\frac{\theta \circ \varphi}{\theta} \in \mathcal{S}(\mathbb{D}).
\]
Here $\es(\D)$ denote the set of all functions that are holomorphic and bounded by one in modulus on $\D$, that is
\[
\es(\D) = \{\psi \in H^\infty(\D): \|\psi\|_\infty: = \sup_{z \in \D} |\psi(z)| \leq 1\}.
\]
The set $\es(\D)$ is known as the \textit{Schur class} and the elements of $\es(\D)$ are called \textit{Schur functions} (see Schur \cite{Schur1, Schur2} and also the monograph \cite{Alpay}).

The proof of the above result, as presented in Section \ref{Mjay1-sec2}, is a simple application of Riesz factorization theorem for $H^2$ functions. Moreover, it is curious to note that several variants of the above result have been used, implicitly, in a number of constructions and proofs in the existing literature (see for instance \cite{Inv2014,Inv2005,Inv2001,Inv2015}). In Section \ref{examples}, we present this point of view by recovering and improving some known results.

In the final section, we point out and correct an error in a corollary of Jones \cite{Inv2005}. On the contrary to the claim of Part 1 of  \cite[Corollary 1]{Inv2005}, in Theorem \ref{thm-error} we prove that for a parabolic automorphism $\vp$ of $\D$, the closed subspace $B_z H^2$ is invariant under $C_{\vp}$, where $B_z$ is the Blaschke product corresponding to the orbit $\{\varphi_{m}(z)\}_{m \geq 0}$ and $z \in \D$ (here $\varphi_m$ denotes the composition of $\vp$ with itself $m$ times).

For general theory of composition operators on $H^2$ we refer the reader to Cowen \cite{CC} and the books by Cowen and MacCluer \cite{Cowen:Book} and Shapiro \cite{Shapiro:Book}.

\section{Invariant subspaces}\label{Mjay1-sec2}

We begin by recalling basic facts about Hardy space and bounded holomorphic functions on $\D$ and refer the reader to Duren \cite[Chapter 2]{Duren:Hpspace} for a more detailed exposition.

Let $\OO(\D)$ denote the set of all holomorphic functions on $\D$. We define the Hardy space $H^2$ as the set of all functions $f \in \OO(\D)$ such that
\[
\|f\|_{2}: =\displaystyle\sup_{0\leq r<1} \Big(\frac{1}{2 \pi} \int_{0}^{2\pi}|f(re^{it})|^{2}\, dt \Big)^{\frac{1}{2}} < \infty.
\]
It is well known (due to Fatou's theorem) that for $f \in H^2$, the \textit{radial limit}
\[
\tilde{f}(e^{it}) := \lim_{r \rightarrow 1^-} f(re^{it}),
\]
exists almost everywhere and $\tilde{f} \in L^2(\partial \D)$ (with respect to the Lebesgue measure on $\partial \D$). In what follows, we will identify $f$ with $\tilde{f}$ and regard $H^2$ as the closed subspace of $L^2(\partial \D)$. Therefore
\[
H^2 = \overline{\mathbb{C}[z]}^{L^2(\partial \D)},
\]
and
\[
\langle f, g\rangle_{H^2} = \frac{1}{2\pi}\int_{0}^{2\pi} f(e^{it}) \overline{g(e^{it})}\,dt \quad \quad (f, g \in H^2).
\]
The space $H^\infty(\D)$ consists of all bounded functions $\psi \in \OO(\D)$. Clearly $H^\infty(\D) \subseteq H^2$, and $H^\infty(\D)$ is a Banach algebra with respect to the uniform norm. Therefore, $\es(\D)$ is the closed unit ball of $H^\infty(\D)$. It is also worth noting that (cf. \cite[Corollary 1.1.24]{Rosenthal:Book})
\[
H^2 \cap L^\infty(\partial \D) = H^\infty(\D).
\]
Recall again that a function $\theta \in \OO(\D)$ is said to be an \textit{inner function} if $|\theta(z)|\leq 1$ for all $z\in \D$ (in particular, $\theta \in H^\infty(\D)$) and its radial limit $|\theta (e^{it})|=1$ a.e. on $\partial \D$. Every inner  function $\theta$ can  be  factored  into a Blaschke product and a singular inner function. That is
\[
\theta = B S,
\]
where the Blaschke product
\[
B(z)= z^m \prod\limits_{n=1}^\infty \frac{|a_n|}{a_n} \frac{a_n-z}{1-\overline{a}_n z} \quad \quad (z\in \mathbb{D}),
\]
for some non-negative integer $m$, is constructed from the zeros of $\theta$ and the singular inner factor
\[
S(z) = c \exp\Big( - \int_0^{2 \pi} \frac{e^{it} + z}{e^{it} - z} \,d\mu(t) \Big) \quad \quad (z \in \D),
\]
for some unimodular constant $c$ and positive measure $\mu$ supported on a  set of Lebesgue measure zero, has no zeros in $\D$. Along the same line, Riesz factorization theorem is enormously useful \cite[Theorem 2.5]{Duren:Hpspace}:

\begin{thm}[Riesz factorization theorem]
\label{factor}
Let $f$ be a non-zero function in $H^2$. Then there exist a Blaschke product $B$ and a function $g \in H^2$ such that $g(z) \neq 0$ for all $z \in \D$ and
\[
f = B g.
\]
Moreover, if $f \in H^\infty(\D)$, then $g \in H^\infty(\D)$ and $\|f\|_\infty=\|g\|_\infty$.
\end{thm}

It is worth noticing that every Blaschke product is an inner function.

Denote by $Z(f)$ the zero set of a holomorphic function $f \in \OO(\D)$. The \textit{multiplicity (or, order)} of $w \in Z(f)$ will be denoted by $mult_f(w)$.

We now return to invariant subspaces of composition operators. \textsf{Throughout this article, $\vp$ will denote a holomorphic self map of $\D$ and $\theta$ will denote an inner function in $H^\infty(\D)$.}

Suppose $C_\vp (\theta H^2) \subseteq \theta H^2$. Then there exists $f \in H^2$ such that
\[
C_{\vp} (\theta 1) = \theta \circ \vp = \theta f.
\]
This yields
\[
Z(\theta) \subseteq Z(\theta \circ \vp),
\]
or equivalently
\[
\vp (Z(\theta)) \subseteq Z(\theta).
\]
More generally, the following easy-to-see remarks adds additional illustration of the concept of zero sets.

\begin{rem}\label{multi-inv}
(1) If $\theta H^2$ is an invariant subspace for $C_\vp$, then
\[
mult_{\theta}(\alpha)\leq mult_{\theta\circ\vp} (\alpha),
\]
for all $\alpha \in Z(\theta)$.

(2) The quotient $\ds\frac{\theta\circ\vp}{\theta}$ defines a holomorphic function on $\mathbb{D}$ if and only if
\[
mult_{\theta}(\alpha)\leq mult_{\theta\circ\vp} (\alpha),
\]
for all $\alpha \in Z(\theta)$.
\end{rem}

The first inequality is merely a necessary condition for $\theta H^2$ to be invariant under $C_{\vp}$ and is not a sufficient condition. A converse of the first remark will be discussed in the next section (see Corollary \ref{cor-2}). Moreover, it is equally evident that the problem of determining effective sufficient conditions, in terms of zero sets of holomorphic functions, is more elusive for zero-free holomorphic functions (like singular inner functions).

Now we are ready to present the central result of this paper.

\begin{thm}\label{main}
The following statements are equivalent:
\begin{enumerate}
 \item[(a)] $\theta H^2$ is an invariant subspace for $C_\vp$.
 \item[(b)] $\ds\frac{\theta\circ\vp}{\theta}\in \es(\D)$.
\end{enumerate}
\end{thm}
\begin{proof}
${\rm (a)}\Rightarrow {\rm (b)}$:  Suppose $\theta H^2$ is an invariant subspace for $C_\vp$. By Remark \ref{multi-inv}, we see that
\[
\ds\frac{\theta\circ\vp}{\theta} \in \OO(\D).
\]
Since $\theta\circ\vp \in\theta H^2$, there exists $f \in H^2$ such that
\[
\theta\circ\vp =\theta f.
\]
It follows that
\[
f =\ds\frac{\theta\circ\vp}{\theta}\in H^2.
\]
Now by Theorem \ref{factor}, there exist a function $g_1 \in H^\infty(\D)$ and a Blaschke product $B_1$ (note that $B_1(z)\equiv 1$ if $Z(\theta)=\emptyset$) such that $g_1(z) \neq 0$ for all $z \in \D$ and
\[
\theta= B_1 g_1.
\]
Since $Z(\theta) \subseteq Z(\theta \circ \vp)$, again by Theorem \ref{factor}, there exist a function $g_2 \in H^\infty(\D)$ and a Blaschke product $B_2$ such that $g_2(z) \neq 0$ for all $z \in \D$ and
\[
\theta\circ\vp = B_1 B_2 g_2.
\]
Since $g_2 \in H^\infty(\D)$ and $\|B_1 B_2\|_{\infty} =1$, as $B_1 B_2$ is an inner function, it follows that
\[
\|g_2\|_\infty = \|B_1 B_2 g_2\|_\infty = \|\theta\circ\vp\|_\infty \leq 1.
\]
Observe
\[
f = \ds\frac{\theta\circ\vp}{\theta}=\frac{B_2 g_2}{g_1} \in H^2.
\]
As $|g_1(e^{it})| = 1$ a.e., by taking the radial limit of both sides, we get
\[
|f(e^{it})|=\left|\frac{g_2(e^{it})B_2(e^{it})}{g_1(e^{it})}\right|=|g_2(e^{it})| ~\mbox{~a.e.}
\]
Hence $f \in H^\infty(\D)$ and $\|f\|_\infty=\|g_2\|_\infty \leq 1$. Therefore $f \in \es(\D)$.

\noindent ${\rm (b)}\Rightarrow {\rm (a)}$: Suppose $\ds\frac{\theta\circ\vp}{\theta}\in \es(\D)$. Then, there exists $f \in \es(\D)$ such that
\[
\theta\circ\vp= \theta f.
\]
Suppose $h \in H^2$. Then
\[
C_\vp(\theta h)= (\theta\circ\vp)\, (h \circ\vp)= \theta f \, (h \circ\vp).
\]
On the other hand,
\[
h \circ\vp \in H^2,
\]
since $C_\vp$ is bounded. As $f \in H^\infty(\D)$, we have $f \, (h \circ\vp)\in H^2$ and hence $C_\vp(\theta h) \in \theta H^2$. This completes the proof of the theorem.
\end{proof}

It is worth noting that the above proof depends on the Riesz factorization theorem on the Hardy space $H^2$. Thus, the above classification result is also valid for $H^p$ spaces on $\D$.

Given the standard factorization of Hardy space functions in a product of a
Blaschke product, a singular inner function, and an outer function, it is clear that
$C_\vp (BH^2) \subseteq BH^2$, if and only if the Blaschke product $B_1$ in the standard factorization
of $B\circ\vp$ is representable as $B_1 = BB_2$, where $B_2$ is a (possibly constant) Blaschke
product, a fact that can be written in terms of multiplicity functions like in the
text of the below corollary.

\begin{cor}\label{cor-2}
Let $B$ be a Blaschke product and let $\vp$ be a holomorphic self map of $\mathbb{D}$. Then the  following statements are equivalent:
\begin{enumerate}
\item $BH^2$ is invariant under $C_\vp$.
\item $mult_B(w) \leq mult_{B\circ\vp}(w)$ for all $w$ in $Z(B)$.
\end{enumerate}
\end{cor}

We refer to Cowen and Wahl \cite[Lemma 8]{Inv2014} for a particular case (where $\vp$ is a non-constant and non-elliptic automorphism) of the above result. Also the special case of inner functions $\vp$ is due to Jones \cite[Lemma 1]{Inv2005}.

Now we proceed to prove a bounded extension problem. Recall that the Hardy space $H^2$ is also a reproducing kernel Hilbert space corresponding to the Szeg\"{o} kernel
\[
K(z, w) = (1 - z \bar{w})^{-1} \quad \quad \quad (z, w \in \D).
\]
For each $w \in \D$, denote by $K(\cdot, w) \in H^2$ the \textit{kernel function} at $w$:
\[
\Big(K(\cdot, w)\Big)(z) = K(z, w) \quad \quad \quad (z \in \D).
\]
The Szeg\"{o} kernel has the following reproducing property:
\[
f(w) = \langle f, K(\cdot, w)\rangle,
\]
for all $f \in H^2$ and $w \in \D$. By using this property, one readily checks that
\[
M_{\psi}^* K(\cdot, w) = \overline{\psi(w)} K(\cdot, w),
\]
and
\[
C_{\vp}^{\ast}K(\cdot, w) = K(\cdot, \vp(w)),
\]
for all $w \in \D$ and $\psi \in H^\infty(\D)$.
\begin{prop}\label{cor}
The following statements are equivalent:
\begin{enumerate}
\item[(a)] $C_{\vp} (\theta H^2) \subseteq \theta H^2$.
\item[(b)] The map
\[
A \left(\overline{\theta(w)}\, K(\cdot, w) \right)= \overline{\theta(\vp(w))}\, K(\cdot, \vp(w)) \quad \quad (w \in \D),
\]
extends to a bounded linear operator on $H^2$.
\item[(c)] $\ds\frac{\theta\circ\vp}{\theta}\in \es(\D)$.
\end{enumerate}
\end{prop}
\begin{proof}
We observe that $C_\vp (\theta H^2) \subseteq \theta H^2$ if and only if
\[
\mbox{ran~} (C_{\vp} M_{\theta}) \subseteq \mbox{ran~} M_{\theta},
\]
which is, by Douglas range inclusion theorem \cite[Theorem 1]{Douglas:Inclusion}, equivalent to
\[
C_\vp M_\theta= M_\theta X,
\]
or equivalently
\[
X^* M_{\theta}^* =M_{\theta}^* C_{\vp}^*,
\]
for some bounded linear operator $X$ on $H^2$. Evaluating each side of the equation by the kernel function $K(\cdot, w)$, $w \in \D$, we get
\[
X^* \left(\overline{\theta(w)}\, K(\cdot, w) \right)= \overline{\theta(\vp(w))}\, K(\cdot, \vp(w)).
\]
Since $\{K(\cdot, w): w \in \D\}$ is a total set in $H^2$, the result follows from Theorem \ref{main}.
\end{proof}

\section{Applications}\label{examples}

We begin by recalling the notion of fixed points of holomorphic self maps. Let $\vp$ be a holomorphic self map of $\D$ and let $w \in \overline{\D}$. We say that $w$ is a \textit{fixed point} \cite[page 50]{Cowen:Book} of $\vp$ if
\[
\lim_{r \rightarrow 1^{-}} \vp(r w) = w.
\]
By a well known result \cite[page 51]{Cowen:Book}, if $w \in \partial \D$ is a fixed point of $\vp$, then
\[
\vp'(w) = \lim_{r \rightarrow 1^{-}} \vp'(r w),
\]
exists as a positive real number or $+\infty$. Now let $\vp$ be an automorphism of $\D$. We say that $\vp$ is:
\begin{enumerate}
\item \textit{elliptic} if it has exactly one fixed point
situated in $\D$,
\item \textit{hyperbolic} if it has two distinct fixed points in $\partial \D$, and
\item \textit{parabolic} if there is only one fixed point in $\partial \D$.
\end{enumerate}

Next we recall the \textit{Denjoy-Wolff theorem}: Let $\vp$ be a holomorphic self map of $\D$. If $\vp$ is not an elliptic automorphism, then there exists $w \in \overline{\D}$ such that $\vp_n$ (the composition of $\vp$ with itself $n$ times) converges to the constant function $w$ uniformly on compact subsets of $\D$. Moreover, $\vp(w) = w$ and (i) $|\vp'(w)| < 1$ if $w \in \D$, and (ii) $0 < \vp'(w) \leq 1$ if $w \in \partial \D$.

The point $w$ is referred to as the \textit{Denjoy-Wolff point} of $\vp$. In connection with the notion of Denjoy-Wolff point and Denjoy-Wolff theorem, we refer the interested reader to \cite[Chapter 2]{Cowen:Book} (also see \cite{D,W}).

By combining Theorem \ref{main} with \cite[Corollary 7]{Inv2014} or \cite[Theorem 2.11]{Inv2015} we obtain the following result concerning shift invariant subspaces generated by atomic singular inner functions.

\begin{thm}\label{Sing.inner}
Let $\vp$ be a holomorphic self map of $\D$ and let $\alpha > 0$. Consider the atomic singular inner function $\theta (z) = e^{\alpha\left(\frac{z+1}{z-1}\right)}$, $z \in \D$. Then the following statements are equivalent:
\begin{enumerate}
\item[(a)]  $\vp(1)= 1$ and $\vp'(1) \leq 1$, that is, $1$ is the Denjoy-Wolff point of $\vp$.
\item[(b)] $\ds\frac{\theta\circ\vp}{\theta} \in \es(\D)$.
\end{enumerate}
\end{thm}

We turn now to a remarkable theorem, due to Matache \cite{Inv2015}, that given a holomorphic self map $\vp$ of $\D$, there exists an inner function $\theta \in H^\infty$ such that $\theta H^2 \subsetneqq H^2$ and
\[
C_{\vp} (\theta H^2) \subseteq \theta H^2.
\]
This is one of the main results of the paper \cite{Inv2015}. Here, we reprove Matache's result. However, our proof is somewhat shorter and simpler.

But before presenting the result, we recall the notion of invariant subspace lattices of operators and make one additional useful observation: For a bounded linear operator $T$ on a Hilbert space $\mathcal{H}$ we denote by $\mbox{Lat}\, T$ the lattice of $T$, that is, the set of all closed invariant subspaces of $T$.

\noindent Now, let $\vp$ is a holomorphic self map of $\D$ and let $a \in \partial \D$. Define $\omega$ and $\psi$, holomorphic self maps of $\D$, by
\[
\omega(z)=\overline{a}z \quad \mbox{and} \quad \psi=\omega\circ \vp \circ \omega^{-1},
\]
for all $z \in \D$. It is easy to see that $a$ is the Denjoy-Wolff point of $\varphi$ if and only if $1$ is the Denjoy-Wolff point of $\psi$. Moreover, if $\theta$ is an inner function, then $\theta H^2 \in \mbox{Lat} \, C_{\psi}$ if and only if (by Theorem \ref{main}) $\theta \circ \psi=  \theta g$ for some $g\in \es(\D)$. On the other hand, $\theta \circ \psi=\theta \circ (\omega\circ \vp \circ \omega^{-1})$ and $g \circ \omega \in \es(\D)$. Hence
\[
(\theta \circ \omega) \circ \vp = \theta g \circ \omega = (\theta \circ \omega) (g \circ \omega),
\]
implies, again by Theorem \ref{main}, that $(\theta \circ \omega) H^2 \in \mbox{Lat} \, C_{\vp}$. In summary, we have the following: (i) $\theta H^2 \in \mbox{Lat} \, C_{\psi}$ if and only if $(\theta \circ \omega) H^2 \in \mbox{Lat} \, C_{\vp}$, and (ii) $a$ is the Denjoy-Wolff point of $\varphi$ if and only if $1$ is the Denjoy-Wolff point of $\psi$.

\begin{thm}\label{thm-nontrivial inv sub}
If $\vp$ is a holomorphic self map of $\D$, then there exists a non-zero closed subspace $\es \subsetneqq H^2$ such that
\[
\es \in \mbox{Lat} \, C_{\vp} \cap \mbox{Lat} \, M_z.
\]
\end{thm}
\begin{proof} Suppose $\vp$ has a fixed point $\alpha$ in $\D$.
Consider the inner function (Blaschke factor)
\[
\theta(z) = \frac{\alpha - z}{1-\overline{\alpha} z} \quad \quad (z \in \D).
\]
Clearly, $\alpha$ is also a zero of $\theta\circ\varphi$ with multiplicity at least one, and so Corollary \ref{cor-2}, we have $C_\varphi (\theta H^2) \subseteq \theta H^2$.

\noindent Finally, suppose $\vp$ does not have any fixed point in $\D$. Then the Denjoy-Wolff point $a$ of $\varphi$ must necessarily lie on $\partial \mathbb{D}$, and so by Theorem \ref{Sing.inner} (along with the remark above), $e^{\alpha\left(\frac{z+a}{z-a}\right)} H^2$ is invariant under $C_\varphi$ for all $\alpha > 0$. This completes the proof of the theorem.
\end{proof}

In the case of elliptic automorphisms of $\mathbb{D}$, Theorem \ref{main} is more definite:

\begin{thm}\label{thm-eigenfunction}
Let $\theta$ be an inner function and $\varphi$ be an elliptic automorphism of $\mathbb{D}$. Then the following statements are equivalent:
\begin{enumerate}
 \item[(a)] $\theta H^2$ is invariant under $C_\varphi$.
 \item[(b)] $\ds\frac{\theta\circ\varphi}{\theta}$ is unimodular constant.
\end{enumerate}
Moreover, in this case, if $w \in \D$ is the unique fixed point of $\varphi$, then
\[
\frac{\theta\circ\varphi}{\theta} \equiv
\begin{cases} \Big(\varphi'(w)\Big)^{\mbox{mult}_{\theta}(w)} & \mbox{if}~ w \in Z(\theta)
\\
\quad 1 & \mbox{otherwise}. \end{cases}
\]
\end{thm}
\begin{proof}
Suppose $\theta H^2$ is invariant under $C_\varphi$. By Theorem \ref{main}, there exists $f\in \es(\D)$ such that $f=\ds\frac{\theta\circ\varphi}{\theta}$. Suppose $w \in \D$ is the unique fixed point of $\varphi$. Define
\[
b_{w}(z) = \frac{w - z}{1 - \bar{w} z} \quad \quad (z \in \D).
\]
Now, if $w \in Z(\theta)$, then there exists an inner function $\theta_1$ such that $\theta_1(w) \neq 0$ and
\[
\theta(z) =
\Big(b_{w}(z) \Big)^{{mult}_{\theta}(w)} \theta_1(z) \quad \quad (z \in \D).
\]
Using this we get
\[
f=\left(\frac{ b_{w} \circ \varphi}{b_{w}}\right)^{{mult}_{\theta}(w)} \frac{\theta_1\circ\varphi}{\theta_1}.
\]
On the other hand
\[
\lim_{z \rightarrow w} \frac{ b_{w} \circ \varphi}{b_{w}} = \varphi'(w),
\]
and $\varphi(w)= w$ implies that
\[
f(w)= \varphi'(w)^{{mult}_{\theta}(w)}.
\]
But, since $\varphi$ is an elliptic automorphism, we have that $|\varphi'(w)| = 1$, and hence $|f(w)|=1$. Then the maximum modulus principle implies that $f\equiv \varphi'(w)^{{mult}_{\theta}(w)}$. Clearly, if $\theta(w)\neq 0$, then $f\equiv f(w)=1$.

\noindent The converse part follows directly from Theorem \ref{main}.
\end{proof}

The above theorem can be reformulated simply as: Let $\theta$ be an inner function and $\varphi$ be an elliptic automorphism of $\mathbb{D}$. Then $C_{\vp} (\theta H^2) \subseteq \theta H^2$ if and only if $\theta$ is an inner eigenfunction of $C_{\vp}$. This result also follows from \cite[Corollary 1.7, Proposition 2.9 and Corollary 2.10]{Inv2015}. However, the present proof is new and somewhat more direct.

The same proof of Theorem \ref{thm-eigenfunction} yields the following result:

\begin{cor}\label{cor-5}
Let $\varphi$  be a holomorphic self map of $\mathbb{D}$ and let $w \in \D$ be the fixed point of $\vp$. Let $\theta$ be an inner function and suppose that $\theta(w)\neq 0$. Then $\theta H^2$ is invariant under $C_\varphi$ if and only if $\theta\circ\varphi=\theta$.
\end{cor}

Now we prove a more definite result on non-automorphic holomorphic self maps.

\begin{cor}
Let $\varphi$  be a non-automorphic and holomorphic self map of $\mathbb{D}$ and let $w \in \D$ be the fixed point of $\vp$. Let $\theta$ be an inner function and suppose that $\theta(w)\neq 0$. Then $\theta H^2$ is invariant under $C_\varphi$ if and only if $\theta$ is an unimodular constant.
In particular, if $\theta$ is a singular inner function, then $\theta H^2$ cannot be invariant under $C_\varphi$.
\end{cor}
\begin{proof}
Suppose $\theta H^2$ is invariant under $C_\varphi$. By Corollary \ref{cor-5}, $\theta\circ\varphi=\theta$, and hence
\[
\theta \circ \varphi_m =\theta,
\]
for all $m \geq 1$ (here $\varphi_m$ denote the composition of $\vp$ with itself $m$ times). Since $\varphi_m$ converges uniformly to the constant function $w$ on every compact subset of $\mathbb{D}$, it follows that $\theta\equiv \theta(w)$. Since $\theta$ is an inner function, we see that $\theta(w)$ is a unimodular constant. The converse part again follows from Theorem \ref{main}.
\end{proof}

\section{Final comments and results}\label{sec-concluding}

We are mainly concerned here with Part 1 of \cite[Corollary 1]{Inv2005}: ``If $\vp$ is a parabolic automorphism then $Lat \, C_{\vp}$ contains no non-trivial $B H^p$.'' This claim is incorrect. Indeed, on the contrary, we prove the following (as always $\varphi_m$ denotes the composition of $\vp$ with itself $m$ times):

\begin{thm}\label{thm-error}
If $\vp$ be a parabolic automorphism of $\D$, then (i) every orbit of $\vp$ is Blaschke summable, and (ii) for each $z\in \D$ we have
\[
B_z H^2 \in \mbox{Lat} \, C_{\vp},
\]
where $B_z$ is the Blaschke product corresponding to the orbit $\{\varphi_{m}(z)\}_{m \geq 0}$.
\end{thm}
\begin{proof}
Let $\vp$ be a parabolic automorphism of $\D$. Suppose
\[
\omega(z)= \frac{1+z}{1-z} \quad \quad (z \in \D).
\]
Then $w$ is a conformal map from $\D$ onto the right half-plane $\mathbb{H}$. Note that
\[
\omega^{-1}(s) = \frac{s-1}{s+1} \quad \quad (s \in \mathbb{H}).
\]
Set
\[
\sigma=\omega\circ \varphi \circ \omega^{-1}.
\]
Then there exists a non-zero real number $b$ such that
\[
\sigma(s)= s + i b \quad \quad (s \in \mathbb{H}),
\]
by the Linear-Fractional Model Theorem (cf. \cite[Section 2.4]{Cowen:Book}). On the other hand,
\[
\varphi_m = \omega^{-1}\circ \sigma_m \circ \omega,
\]
for all $m$, and hence
\[
\begin{split}
1-|\varphi_m(z)|^2 & = 1-|\omega^{-1}(\sigma_m (\omega(z)))|^2
\\
& = 1 - \left|\frac{\sigma_m(\omega(z)) - 1}{\sigma_{m}(\omega(z)) + 1}\right|^2
\\
& = \frac{4 \mbox{~Re~} \Big( \sigma_m(\omega(z))\Big)}{|\sigma_m(\omega(z)) + 1|^2},
\end{split}
\]
for all $z \in \D$. Now we fix $z \in \D$ and let $\omega(z)=u+i\, v$. Then
\[
\sigma_m(\omega(z)) = \omega(z) + i m b = u + i(mb+v),
\]
for all $m$. It follows that
\[
|\sigma_m(\omega(z))+1|^2 = (1+u)^2+(mb+v)^2,
\]
and hence
\[
1-|\varphi_m(z)|^2= \frac{4u}{(mb+v)^2+(1+u)^2} \sim \frac{4u}{b^2m^2},
\]
for large $m$. Therefore
\[
\sum_{m} 1-|\varphi_m(z)|^2 < \infty \quad \quad (z\in \D).
\]
Hence $|\varphi_{m}(z)| \geq |\varphi_m(z)|^2$ for all $m$ yields that
\[
\sum_{m} 1-|\varphi_m(z)| < \infty,
\]
that is, the orbit $\{\vp_m(z)\}_{m \geq 0}$ of $\vp$ at $z\in \D$ is Blaschke summable. The second part follows from the first and Corollary \ref{cor-2}. This completes the proof of the theorem.
\end{proof}

From the above proof it is now evident that the estimate
\[
``1-|\varphi_n(z)|^2 \sim \frac{c}{n}",
\]
in the proof of Part 2 of \cite[Lemma 3]{Inv2005} is incorrect.

To conclude, we remark that a Schur function always admits a fractional linear transformation representations in the following sense: Given $\vp \in \es(\D)$, there exist a Hilbert space $\mathcal{H}$ and a unitary (/isometry/co-isometry/contractive) matrix
\[
U = \begin{bmatrix} a & B \\ C & D \end{bmatrix} : \mathbb{C} \oplus \mathcal{H} \rightarrow \mathbb{C} \oplus \mathcal{H},
\]
such that
\[
\vp(z) = a + z B(I - z D)^{-1} C \quad \quad (z \in \D).
\]
This point of view has proved extremely fruitful in understanding the structure of composition operators (cf. \cite{MJ}). In the context of Theorem \ref{main}, a number of questions arise naturally at this point. For instance, a natural question arises as to whether one can relate the fractional linear transformations of $\vp$ and $\theta$ with the fractional linear transformation of $\frac{\theta\circ\vp}{\theta}$. We hope to return to this theme in future work.

\vspace{0.3in}

\NI\textit{Acknowledgement:} The authors are very grateful to the referee for careful reading of the paper and valuable suggestions and comments. The authors would also like to thank Sushil Gorai for enlightening comments and suggestions. The first and second authors thanks the National Board for Higher Mathematics (NBHM), India, for providing financial support  to carry out this research. The research of the third named author is supported in part by National Board of Higher Mathematics (NBHM) grant NBHM/R.P.64/2014, and the Mathematical Research Impact Centric Support (MATRICS) grant, File No: MTR/2017/000522 and Core Research Grant, File No: CRG/2019/000908, by the Science and Engineering Research Board (SERB), Department of Science \& Technology (DST), Government of India.

\end{document}